\documentclass[12pt,a4paper]{amsart}
\usepackage{amsmath,amssymb,amsfonts,amsthm,enumerate}

\newtheorem{lemma}{Lemma}
\newtheorem{remark}{Remark}

\theoremstyle{remark}\newtheorem{as}{Assumption}
\theoremstyle{plain}\newtheorem{theorem}{Theorem}

\newcommand{\N}{\mathbb{N}}

\newcommand{\biz}{\textbf{Proof.\ }}
\newcommand{\ri}{\right)}
\newcommand{\lef}{\left(}

\newcommand{\E}{\mathbb{E}}

\newcommand{\lk}{\left\lbrace}
\newcommand{\rk}{\right\rbrace}

\newcommand{\Prob}{\mathbb{P}}

\newcommand{\lsz}{\left[}
\newcommand{\rsz}{\right]}

\newcommand{\F}{\mathcal{F}}
\newcommand{\tv}{\rightarrow \infty}

\newcommand{\al}{\alpha}
\newcommand{\pls}{\!+\!}
\newcommand{\mns}{\!-\!}
\newcommand{\var}{\mathop{\textrm{Var}}}

\newcommand{\veps}{\varepsilon}
\newcommand{\vph}{\varphi}

\newcommand{\kesz}{\hfill $\square$}
\newcommand{\Cal}{\mathcal}

\title{A random model of publication activity}
 
\author{\'Agnes Backhausz}
\address{Department of Probability Theory and Statistics\\
E\"otv\"os Lor\'and University\\P\'azm\'any P.~s.\ 1/C, H-1117 Budapest, Hungary} 
 \email{agnes@cs.elte.hu}
 \author{Tam\'as F.~M\'ori}
 \address{Department of Probability Theory and Statistics\\
 E\"otv\"os Lor\'and University\\P\'azm\'any P.~s. 1/C, 
 H-1117 Budapest, Hungary}
 \email{moritamas@ludens.elte.hu}
 \dedicatory{\upshape Department of Probability Theory and Statistics,
\\ E\"otv\"os Lor\'and University\\ P\'azm\'any P.~s. 1/C,
 H-1117 Budapest, Hungary\\                         
 \textit{E-mail address:} \texttt{agnes@cs.elte.hu, moritamas@ludens.elte.hu} }

 

\thanks{The European Union and the European Social Fund have provided
financial support to the project under the grant agreement no. T\'AMOP
4.2.1./B-09/KMR-2010-0003.}
  
\date{7 November 2014}
\begin{document}

\begin{abstract}
We examine a random structure consisting of objects with 
positive weights and evolving in discrete time steps. It generalizes 
certain random graph models. We prove almost sure convergence for the 
weight distribution and show scale-free asymptotic behaviour. Martingale 
theory and renewal-like equations are used in the proofs. 
\end{abstract}
\maketitle
\thispagestyle{empty}

\section{Introduction}

In this paper we examine a dynamic model inspired by scientific 
publication activity and networks of coauthors. However, the 
model contains many simplifying assumptions that are not valid 
in reality. We still use the terminology of publications for 
sake of simplicity.

The model consists of a sequence of researchers. Each of them has a
positive weight which is increasing in discrete time steps. The
weights reflect the number and importance of the researcher's
publications. One can think of cumulative impact factor for instance.

We start with a single researcher having a random positive weight. 
At the $n$th step a new publication is born. The number of its authors  
is randomly chosen. Then we select the authors, that is, one of the
groups of that size; the probability that a given group is chosen is
proportional to the sum of the weights of its members. After that the
weights of the authors of the new publication are increased by random
bonuses. Finally, a new researcher is added to the system with a
random initial weight.  

This is a preferential attachment model; one can see that authors with 
higher weights have larger chance to be chosen and increase their
weights when the new publication is born. 

We are interested in the weight distribution of the model. That is,
for fixed $t>0$,  we consider the ratio of authors of weight larger
than $t$, and study the asymptotic behaviour of this quantity as the
number of steps goes to infinity.  

Our main results (Section 3) include the almost sure convergence of
the ratio of authors of weight larger than $t$ under suitable
conditions; first, when all weights are integer valued, then assuming
that these random variables have continuous distribution. In both
cases we describe the limiting sequence or function and determine its
asymptotics. They are polynomially decaying under suitable conditions, 
thus our model shows scale-free behaviour.  

The proofs of the almost sure convergence are based on the methods of
martingale theory, while the polynomial decay of the asymptotic weight
distribution follows from the results of \cite{rek} about renewal-like 
equations. See Section 4 for the details.

This model generalizes some random graph models. To see this, assume
that every publication has only one author, and at each step, when a
publication is born, connect its author to the new one with an
edge. We get a random tree evolving in time. 

In the particular case where the initial weights and author's bonuses 
are always equal to $1$, we get the Albert--Barab\'asi random tree
\cite{ba}. The neighbour of the new vertex is chosen with
probabilities proportional to the degrees of the old
vertices. Similarly, if the initial weights and the bonuses are fixed,
but they are not necessarily equal to each other, we get random trees
with linear weights \cite{pittel}, sometimes called generalized plane
oriented recursive trees. In these cases the asymptotic degree
distribution is well-known.

\section{Notations and assumptions}
\subsection{Notations}

Let the label of the only researcher being present in the beginning be
$0$; the label of the researcher coming in the $n$th step is $n$.

$X_i$ is the initial weight of researcher $i$ for $i=0, 1,\ldots$. We
suppose that $X_0, X_1,\ldots $ are independent, identically
distributed positive random variables.    

$\nu_n$ is the number of coauthors at step $n$. This is an integer
valued random variable for each $n$. Obviously $\nu_n\leq n$ must hold
for all $n\geq 1$. On the other hand, for technical reasons we also
assume that $\nu_n\geq 1$ for all $n\geq 1$. Since the authors'
weights are not necessarily increased, this may be supposed without
loss of generality.  

Given that $\nu_n=k$, a group of size $k$ is chosen randomly from
researchers $0, \ldots, n-1$. The probability that a given group is
chosen is proportional to the total weight of the group. The selected 
researchers will be the authors of the $n$th paper.  

Let $Y_{n,1}, Y_{n,2}, \ldots, Y_{n, \nu_n}$ be nonnegative random
variables. These are the authors' bonuses at step $n$. That is, the
weight of the $i$th coauthor of the $n$th paper is increased by 
$Y_{n,i}$. The order of the coauthors is the natural order of the
labels.   

Let $Z_n$ be the total weight of the $n$th paper; that is, 
$Z_n=Y_{n,1}+Y_{n,2}+ \ldots+ Y_{n,\nu_n}$ for $n\geq 1$. 

$W\lef n, i\ri$ denotes the weight of author $i$ after step $n$ for
$i=0, \ldots, n$. This is equal to $X_i$ plus the sum of all bonuses
$Y_{j,\ell}$ for which author $i$ is the $\ell$th author of the $j$th
paper $\lef\ell=1,\ldots, \nu_j, \ j=1,2,\ldots,n\ri$.  

Let $S_n$ be the total weight after $n$ steps; namely, 
\[
S_n=W\lef n,0\ri+\ldots +W\lef n,n\ri=X_0+\ldots+X_n+Z_1+\ldots+Z_n.
\]

$X,\ \nu,\ Y_n,\ Y$ and $Z$ are random variables. $X$ is equal to
$X_0$ in distribution, and $Y_n$ is equal to $Y_{n,1}$ in distribution
for $n\geq 1$. The other random variables will be determined later by
the assumptions. 

Finally, $\F_n$ is the $\sigma$-algebra generated by the first $n$
steps; $\F_n^+=\sigma\lk \F_n, \nu_{n+1}\rk$. 

Throughout this paper $\mathbb I(A)$ denotes the indicator of event
$A$. We say that two sequences $\lef a_n\ri$, $\lef b_n\ri$ are
asymptotically equal $\lef a_n\sim b_n\ri$, if they are positive
except finitely many terms, and $a_n/b_n\to 1$ as $n\tv$. A sequence
$\lef a_n\ri$ is exponentially small if $\left\vert a_n\right\vert\leq
q^n$ holds for all sufficiently large $n\in \N$ for some $0<q<1$. 

\subsection{Assumptions}

Now we list the assumptions on the model. 

\begin{as} \label{puba1}$X_0, X_1, \ldots$ are independent,
identically distributed. The initial weights $X_n$, and the triplets
$\left(\F_{n-1},\left(Y_{n,1},\ldots,Y_{n,\nu_n}\right),\nu_n\right)$ 
are independent $\lef n=1, 2,\ldots\ri$.
\end{as}

\begin{as} \label{puba2}$X$ has finite moment generating function.
\end{as}

\begin{as}\label{puba3}
$\nu_n$ and $\left(Y_{n,1},\ldots,Y_{n,\nu_n}\right)$ are independent
of $\F_{n-1}$ for $n\geq 1$. 
\end{as}

\begin{as}\label{puba4}
$\nu_n\rightarrow \nu$ in distribution as $n\tv$; in addition,
$\E\nu_n\rightarrow\E\nu <\infty$ and $\E\nu_n^2\rightarrow\E\nu^2 
<\infty$ hold. 
\end{as}

Recall that $\nu_n\leq n$. Assumption 4 trivially holds if
$\nu$ is a fixed random variable with finite second moment, and the
distribution of $\nu_n$ is identical to the distribution of $\min\lef
n,\nu\ri$, or to the conditional distribution of $\nu$ with respect to
$\lk\nu\leq n\rk$.  

\begin{as}\label{puba5}
The conditional distribution of $\left(Y_{n,1},\ldots,Y_{n,\nu_n}
\right)$, given $\nu_n=k$, does not depend on $n$. Moreover, the 
components are conditionally interchangeable, given $\nu_n=k$. 
\end{as}

\begin{as}\label{puba6}
$Z_n$ has finite expectation.
\end{as}

Now we know that $\lef \nu_n,\ Y_n,\ Z_n\ri\rightarrow \lef \nu,\ Y,\
Z\ri$ in distribution as $n\tv$, where $Y$ and $Z$ are random
variables. We need that they also have finite moment generating
functions, and they are not degenerate. 

\begin{as}\label{puba7}
 $Y$ and $Z$ have finite moment generating functions.
\end{as}

\begin{as}\label{puba8}
$X_n$, $Y_n$, $X$ and $Y$ are positive with positive probabilities for
every $n=1,2,\dots$ . In addition, if $Y$ is integer-valued, then the
greatest common divisor of the set $\lk i: \Prob\lef Y=i\ri>0\rk$ is
equal to $1$.  
\end{as}

The condition on the positivity of $X_n$ and $Y_n$ is not
crucial. The positivity of $X$ and $Y$ implies that the same holds for
$X_n$ and $Y_n$ if $n$ is large enough; we may assume this for all $n$
without loss of generality. On the other hand, if $\lef Y_n\ri$ is
identically equal to $0$, that is, there are no bonuses at all, then
the model only consists of the sequence of independent and identically
distributed initial weights $X_n$, and the problem of empirical weight
distribution becomes trivial. The last part of this assumption
excludes periodicity.

There are two important particular cases satisfying all of our
conditions.  In the first one the weight of the paper is equally
distributed among the authors. That is, $Z_1, Z_2, \ldots$ are
independent identically distributed random variables, and
$Y_{n,1}=\ldots=Y_{n,\nu_n}=Z_n/\nu_n$. The other option is that every
author gets the total bonus, regardless the number of coauthors. More
precisely, $Y_1,Y_2,\ldots$ are independent and identically
distributed, and $Y_{n,1}=\ldots=Y_{n,\nu_n}=Y_n$, thus $Z_n=\nu_nY_n$. 

\section{Main results}

\subsection*{Discrete weight distribution}
Suppose first that $X, Y_1, Y_2, \ldots$ are nonnegative integer valued 
random variables. Let $\xi_n\lef j \ri$ denote the number of researchers 
of weight $j$ after $n$ steps, that is, 
\[
\xi_n\lef j\ri=\bigl\vert \lk 0\leq i\leq n: W\lef n,i\ri=
j\rk\bigr\vert,\quad j, n=1,2,\ldots\,.
\] 

The first theorem is about the almost sure behaviour of this quantity.  
 
\begin{theorem}\label{pubt1}
$\dfrac{\xi_n(j)}{n}\rightarrow x_j$ almost surely as $n\tv$ with
positive constants $x_j$, $j=1,2,\ldots$. The sequence $\lef x_j\ri$
satisfies the recursion 
\begin{equation}\label{pubrec}
x_j=\frac{\sum\limits_{i=1}^{j-1} x_{j-i}\biggl[\dfrac{(j\mns i) 
\Prob(Y=i)}{\E X\pls\E Z}+\E\bigl((\nu\mns 1)\mathbb I(Y=i)\bigr)
\biggr]+\Prob(X=j)}{\alpha j+\beta+1}\, ,
\end{equation}               
where $\alpha=\dfrac{\Prob(Y>0)}{\E X\pls\E Z}$,\ \ 
$\beta=\E\bigl((\nu\mns 1)\mathbb I(Y>0)\bigr)$.
\end{theorem}

The second theorem describes the asymptotic behaviour of the sequence
$\lef x_j\ri$. 
\begin{theorem}\label{pubt2}
We have $x_j\sim C\,j^{-\gamma}$ as $j\tv$, where $C$ is a
positive constant, and 
\[
\gamma=\frac{\E X\pls\E Z}{\E Y}\pls 1.
\]
\end{theorem}

\subsection*{Continuous weight distribution}

Now we assume that the distribution of $X$ and the conditional
distributions of $Y_n\mid\nu_n=k$ are continuous for $k=1,2,\dots, n,\ 
n=1,2,\ldots$\,. This implies that the distribution of $Y_n$ is
continuous. Moreover, since the conditional distribution does not
depend on $n$ according to Assumption \ref{puba5}, the distribution of
$Y$ is also continuous.  

Let $F(t)=\Prob(Y>t)$, $H(t)=\E\bigl((\nu-1)\mathbb I(Y>t)\bigr)$, and
\[
L(t,s)=\frac{sF(s)+t(1-F(s))}{\E X+\E Z}-H(s),\quad 0\le s\le t. 
\]
It is clear that $L(t,s)$ is continuous, and, being the difference of
two increasing functions, it is of bounded variation for fixed $t$. 

This time $\xi_n(t)$ denotes the number of researchers with weight more 
than $t$ after $n$ steps.  
\[
\xi_n\lef t\ri=\bigl\vert \lk 0\leq i\leq n: W\lef n,
  i\ri>t\rk\bigr\vert, \quad t>0,\ n=1,2,\ldots\,.
\]

\begin{theorem}\label{pubt3}
$\dfrac{\xi_n(t)}{n}\to G(t)$ almost surely, as $n\tv$, where $G(t)$
is the solution of the following integral equation.
\begin{equation}\label{pube0}
G(t)=\frac{{\displaystyle\int_0^t} G(t-s)\,d_sL(t,s)+H(t)+\Prob(X>t)}
{\dfrac{t}{\E X+\E Z}+\E \nu}
\end{equation}
for $t>0$, and $G(0)=1$.
\end{theorem}

Adding some extra conditions we can obtain results on the asymptotic
behaviour of $G$. 

\begin{theorem}\label{pubt4}
Suppose that the distribution of $Y$ is absolutely continuous. Then we
have $G(t)\sim C\,t^{-\gamma}$ as $t\tv$, where $C$ is a
positive constant, and  
\[
\gamma=\frac{\E X+\E Z}{\E Y}.
\] 
\end{theorem}

\begin{remark}
The difference of the exponents in the discrete and continuous cases
is due to the difference in the definitions. Namely, in the first case
$\xi_n$ denotes the weight distribution, while in the second case it
stands for the complementary cumulative weight distribution function.  
\end{remark}

\section{Proofs}

First we prove some propositions we will often use in the sequel.  

\begin{lemma}\label{publ1}
Let $(\mathcal F_n)$ be a filtration, $(\xi_n)$ a nonnegative adapted
process. 
Let $(w_n)$ be a regularly varying sequence of positive numbers with
exponent $\mu>-1$. Suppose that 
\begin{equation}\label{lemmafelt}
\mathbb \E\bigl((\xi_n-\xi_{n-1})^2\bigm|\mathcal F_{n-1}\bigr)
=O\left( n^{1-\delta+2\mu}\right)
\end{equation}
holds with some $\delta>0$. Let $(u_n)$, $(v_n)$ be nonnegative
predictable processes such that $u_n<n$ for all $n\geq 1$. 

$(a)$ Suppose that 
\[
\mathbb \E(\xi_n\mid\mathcal F_{n-1})\le
\Bigl(1-\dfrac{u_n}{n}\Bigr)\xi_{n-1}+v_n,
\]
and $\lim_{n\tv} u_n=u$, $\limsup_{n\tv} v_n/w_n\le v$ with 
some random variables $u>0,\ v\geq 0$. Then 
\begin{equation*}
\limsup_{n\rightarrow\infty}\frac{\xi_n}{nw_n}\le \frac{v}{u+\mu+1}
\quad a.s.
\end{equation*} 

$(b)$ Suppose that 
\[
\mathbb \E(\xi_n\mid\mathcal F_{n-1})\ge
\Bigl(1-\dfrac{u_n}{n}\Bigr)\xi_{n-1}+v_n,
\]
and $\lim_{n\tv}u_n=u$, $\liminf_{n\tv} v_n/w_n\ge v$ with some random
variables $u>0,\ v\geq 0$. Then  
\begin{equation*}
\liminf_{n\rightarrow\infty}\frac{\xi_n}{nw_n}\ge \frac{v}{u+\mu+1}
\quad a.s.
\end{equation*}
\end{lemma}

This is a stochastic counterpart of a lemma of Chung and Lu
\cite{chung}. We will often apply this proposition with the sequence
$w_n\equiv 1$  and $\mu=0$. 

\biz Suppose first that $v$ is strictly positive. Let $\mathcal F_0$
be the trivial $\sigma$-algebra, $\xi_0=0$, and
\[
c_n=\prod_{i=1}^n\Bigl(1-\frac{u_i}{i}\Bigr)^{-1},\quad n\geq
1.
\]

We have
\[\log c_n=\sum_{i=1}^n\frac{u_i}{i}\bigl(1+o(1)\bigr)=
u\sum_{i=1}^n\frac{1+o(1)}{i}\,.
\]
Hence for all $t>1$ we get that $\lim_{n\tv}
(\log c_{[tn]}-\log c_n)=u\log t$. That is,  $\lef c_n\ri$ is regularly 
varying with exponent $u$. It is clear that
\begin{equation}\label{pube1}\E\bigl(c_n\xi_n\bigm|\mathcal F_{n-1}\bigr)\leq
c_{n-1}\xi_{n-1}+c_nv_n.\end{equation}
Therefore $c_n\xi_n$ is a submartingale. Consider the Doob
decomposition $c_n\xi_n=M_n+A_n$, where
\[
M_n=\sum_{i=1}^n\left(c_i\xi_i-
\E\bigl(c_i\xi_i\bigm|\mathcal F_{i-1}\bigr)\right)
\]
is a martingale, and
\[
A_n=\sum_{i=1}^n\left(\E\bigl(c_i\xi_i\bigm|\mathcal F_{i-1}\bigr)-
c_{i-1}\xi_{i-1}\right).
\]
From inequality \eqref{pube1} it follows that
\[
A_n\leq\sum_{i=1}^n c_iv_i.
\]

Consider the increasing process in the Doob decomposition of the
square of the martingale $(M_n)$. Using condition \eqref{lemmafelt} we
get that 
\begin{align*}
B_n&=\sum_{i=1}^n\var\bigl(c_i\xi_i\bigm|\mathcal F_{i-1}\bigr)
=\sum_{i=1}^n\var\bigl(c_i(\xi_i-\xi_{i-1})\bigm|\mathcal F_{i-1}
\bigr)\\ 
&\leq\sum_{i=1}^n c_i^2\,\E\bigl((\xi_i-\xi_{i-1})^2\bigm|\mathcal
F_{i-1}\bigr)=O\Biggl(\;\sum_{i=1}^n i^{1-\delta+2\mu}c_i^2\Biggr).
\end{align*}

Since $n^{1-\delta+2\mu}c_n^2$ is still regularly varying with exponent
$2u+1-\delta+2\mu$, it follows that $B_n=O\bigl(n^{2-\delta+2\mu}c_n^2\bigr)$
(see e.g. \cite{bingham,bojanic}). Hence, by 
Propositions VII-2-3 and VII-2-4 of \cite{[Ne75]}, we have
\[
M_n=O(B_n^{1/2+\veps}\bigr)=O\bigl(n^{(2-\delta+2\mu)(1/2+\veps)}
c_n^{1+2\veps}\bigr)=o\bigl(nc_nw_n\bigr)\quad a.s.,
\]
for all $0<\veps<\dfrac{\delta}{4(u+1+\mu)}$\,.

On the other hand, using the fact $u+\mu>-1$, and the results of
\cite{bingham, bojanic} on regularly varying sequences  
we obtain that
\[
A_n\leq\sum_{i=1}^nc_iv_i\le\bigl(1+o(1)\bigr)\,v\sum_{i=1}^nc_iw_i\sim
v\,\frac{nc_nw_n}{u+\mu+1}
\]
almost surely, as $n\tv$. This implies that 
\[
c_n\xi_n\le\bigl(1+o(1)\bigr)\frac{v}{u+\mu+1}\,nc_n w_n,
\]
thus the proof of part $(a)$ is complete for positive $v$.

The general case of nonnegative $v$ can be deduced from the positive
case by noticing that 
\[
\mathbb \E(\xi_n\mid\mathcal F_{n-1})\le
\Bigl(1-\dfrac{u_n}{n}\Bigr)\xi_{n-1}+\max\lef v_n, \veps\ri
\] 
for arbitrary $\veps>0$.

The proof of part (b) is similar. In this case
\[
A_n\ge\sum_{i=1}^nc_iv_i\sim\frac{v}{u+\mu+1}\,nc_nw_n,
\]
a.s. on the event $\{v>0\}$. Hence, using $c_n\xi_n\sim A_n$, we get
that 
\[
c_n\xi_n\ge\frac{v}{u+\mu+1}\,nc_nw_n\bigl(1+o(1)\bigr).
\] 
On the event $\lk v=0\rk$ the inequality trivially holds. \qed

\begin{lemma}\label{publ2}
The conditional probability that an author of weight $j$ is chosen,
given $\F_n^+$ and $\nu_{n+1}=k$, is equal to 
\[
\frac{k-1}{n}+\frac{n+1-k}{n}\cdot\frac{j}{S_{n}}
=\frac{k-1}{n}\lef1-\frac{j}{S_n}\ri+\frac{j}{S_n}.
\] 
\end{lemma}

\biz Consider those groups of size $k\ge 2$ that contain researcher $i$
($0\leq i\leq n$). There are $\binom{n}{k-1}$ of them, because the
total number of researchers is $n+1$. Researcher $i$ belongs to all of
them, while the other researchers belong to $\binom{n-1}{k-2}$ of
those groups. Therefore the total weight of these groups can be
obtained in the following way.  

\begin{align*}
\sum_{\substack{H\subset\lk 0,\ldots,n\rk\\|H|=k,\,i\in H}}
\ \sum_{j\in H}W(n,j)&=\binom{n}{k-1}W(n,i)+
\sum_{j\neq i}\binom{n-1}{k-2}W(n,j)\\
&=\binom{n-1}{k-1}W(n,i)+\binom{n-1}{k-2}S_n.
\end{align*}

On the other hand, the total weight of all groups of size $k$ is
given by
\[
\binom{n}{k-1}S_n.
\]

Hence the conditional probability that researcher $i$ participates in
the $\lef n+1\ri$st paper given that it has $k$ authors is equal to 
\[
\frac{k-1}{n}+\frac{n-k+1}{n}\cdot\frac{W(n,i)}{S_n}
=\frac{k-1}{n}\lef 1-\frac{W(n,i)}{S_n}\ri+\frac{W(n,i)}{S_n}.
\]
This obviously holds for $k=1$ as well. 
\qed

\subsection*{Proof of Theorem \ref{pubt1}.}

Recall that in Theorem \ref{pubt1} we assumed that $X, Y_1, Y_2, \ldots$
are integer valued random variables. Let us introduce
\[
H\lef i \ri=\E\lef \lef \nu-1\ri \mathbb I\lef Y=i\ri\ri,
\]
then $\beta=\sum_{i=1}^{\infty} H(i)$.

We prove the theorem by induction on $j$. The following argument is
valid for all $j=1, 2, \ldots$.  
For $j>1$ we will use the induction hypothesis.  

At each step the number of authors of weight $j$ may change due to 
the following events. 
\begin{itemize}
\item
A given author of weight $j$ is chosen and he gets positive bonus.
\item
A given author of weight $j-i$ is chosen and his bonus is equal to
$i$.  
\item
The initial weight of the new author is $j$.
\end{itemize}
Therefore Lemma \ref{publ2} implies that
\begin{multline}\label{pube1.1}
\E\bigl(\xi_n(j)\bigm|\mathcal F_{n-1}^+\bigr)\\
=\xi_{n-1}(j)\biggl[1-\Prob\bigl(Y_n>0\bigm|\mathcal F_{n-1}^+\bigr)
\Bigl(\frac{\nu_n-1}{n-1} + \frac{n-\nu_n}{n-1}\cdot
\frac{j}{S_{n-1}}\Bigr)\biggr]\\
+\sum_{i=1}^{j-1}\xi_{n-1}(j-i)\Prob\bigl(Y_n=i\bigm|
\mathcal F_{n-1}^+\bigr)\Bigl(\frac{\nu_n-1}{n-1}+
\frac{n-\nu_n}{n-1}\cdot\frac{j-i}{S_{n-1}}\Bigr)\\
+\Prob(X_n=j).
\end{multline}
Recall that $\nu_n\geq 1$ is assumed.

We introduce the time-dependent versions of the already defined 
quantities. Namely, 
\[H_n(i)=\E\bigl((\nu_n-1)\mathbb I(Y_n=i)\bigr); \
\beta_n=\sum_{i=1}^n H_n(i)=\E\bigl((\nu_n-1)\mathbb I(Y_n>0)\bigr).
\]

Let us take conditional expectation given $\mathcal F_{n-1}$ in both
sides of \eqref{pube1.1}. Then we get that   
\begin{multline}\label{pube1.2}
\E\bigl(\xi_n(j)\bigm|\mathcal F_{n-1}\bigr)\\
=\xi_{n-1}(j)\biggl[1-\frac{\beta_n}{n-1}-\Bigl(\Prob(Y_n>0)-
\frac{\beta_n}{n-1}\Bigr)\frac{j}{S_{n-1}}\biggr]\\
+\sum_{i=1}^{j-1}\xi_{n-1}(j-i)\biggl[\frac{H_n(i)}{n-1}+
\Bigl(\Prob(Y_n=i)-\frac{H_n(i)}{n-1}\Bigr)\frac{j-i}{S_{n-1}}\biggr]\\
+\Prob(X_n=j) \quad (j,n=1,2,\ldots).
\end{multline}

We are going to apply Lemma \ref{publ1} to the sequence $\lef
\xi_n(j)\ri$ with $w_n\equiv 1$ and $\mu=0$. It is clear that  
$|\xi_n(j)-\xi_{n-1}(j)|\le \nu_n+1$, hence 
\[
\E\bigl((\xi_n(j)-\xi_{n-1}(j))^2\bigm|\mathcal
F_{n-1}\bigr)\le \E(\nu_n+1)^2=O(1).
\]
Thus, condition \eqref{lemmafelt} on the differences of the sequence
$\xi_n(j)$ is satisfied.
Moreover, as $n\tv$, we have 
\[
u_n=n\biggl[\frac{\beta_n}{n-1}+\Bigl(\Prob(Y_n>0)-\frac{\beta_n}{n-1}\Bigr)
\frac{j}{S_{n-1}}\biggr]\to \beta+\alpha j.
\]
Note that $\al>0$ because of Assumption \ref{puba8}.

Though the random variables $Z_1, Z_2,\ldots$ are not necessarily
identically distributed, they satisfy the following conditions.
\[
\sum_{n=1}^{\infty}\frac{\var\lef Z_n\ri}{n^2}<\infty,\quad 
\lim_{n\tv} \frac1n\sum_{i=1}^n \E Z_i=\E Z.
\]

Therefore Kolmogorov's theorem (Theorem 6.7. in \cite{petrov}) can be
applied. We get that $S_n\sim n\lef \E X+\E Z\ri$ almost surely as $n\tv$. 
Using this, and also the induction hypothesis when $j>1$, we conclude
that   
\begin{multline*}
v_n=\sum_{i=1}^{j-1}\xi_{n-1}(j-i)\biggl[\frac{H_n(i)}{n-1}+
\Bigl(\Prob(Y_n=i)-\frac{H_n(i)}{n-1}\Bigr)\frac{j-i}{S_{n-1}}\biggr]
+\Prob(X_n=j)\\
\to \sum_{i=1}^{j-1}x_{j-i}\biggl[H(i)+\Prob(Y=i)\,
\frac{j-i}{\E X+\E Z}\biggr]+\Prob(X=j),
\end{multline*}
as $n\to\infty$. 

From equations \eqref{pube1.1} and \eqref{pube1.2} one can see that
$\lef u_n\ri$ and $\lef v_n\ri$ are nonnegative predictable processes.
Moreover, $u_n<n$ if $n$ is large enough, because then $\nu_n<n$ and
$j<S_{n-1}$. We have also seen that the limit of
$\lef u_n\ri$ is positive. Hence, by Lemma \ref{publ1}, the induction
step and the proof of Theorem \ref{pubt1} is complete.     
\kesz

\subsection*{Proof of Theorem \ref{pubt2}.}

Write recursion \eqref{pubrec} in the following form.
\[
x_j=\sum_{i=1}^{j-1}w_{j,i}x_{j-i}+r_j,
\]
where for $i,j\geq 1$ we set
\[
w_{j,i}=\frac{\biggl[\dfrac{(j\mns i)\mathbb \Prob(Y=i)}
{\mathbb \E X\pls\mathbb \E Z}+\mathbb \E\bigl((\nu\mns 1)
\mathbb I(Y=i)\bigr)\biggr]}{\alpha j+\beta+1}\,,
\]
and 
\[ 
r_j=\frac{\Prob\lef X=j\ri}{\alpha j+\beta+1}\,.
\]

In order to apply Theorem 1 of \cite{rek} we try to find sequences 
$\lef a_i\ri$, $\lef b_i\ri$, $\lef c_{j,\,i}\ri$ such 
that $w_{j,i}=a_i\pls \frac{b_i}{j}\pls c_{j,\,i}$ holds, then we have
to check that $a_i$, $b_i$, $c_{i,j}$, $r_i$ satisfy the 
following conditions.
\begin{enumerate}[(i)]
\item
$a_i\ge 0$ for $i\geq 1$, and the greatest common divisor of the
set $\lk i:  a_i>0\rk$ is  $1$;  
\item
$r_i$ is nonnegative, and not identically zero;
\item\label{iii}
there exists $z>0$ such that
\begin{gather*}
1<\sum_{i=1}^\infty a_iz^i<\infty,\qquad
\sum_{i=1}^\infty|b_i|z^i<\infty,\\
\sum_{i=1}^\infty\sum_{j=1}^{i-1}|c_{i,j}|z^j<\infty,\qquad
\sum_{i=1}^\infty r_iz^i<\infty.
\end{gather*}
\end{enumerate}

Therefore we set
\begin{equation*}
a_i=\lim_{j\to\infty}w_{j,i}=\frac{\Prob(Y=i)}{\alpha(\E X+\E Z)}=
\Prob(Y=i\mid Y>0), \quad i=1, 2,\dots\,,
\end{equation*}
then we define
\[
b_i=\lim_{j\to\infty}j(w_{j,i}-a_i)=\frac{1}{\alpha}
\Bigl[H(i)-(\alpha i+\beta+1)a_i\Bigr].
\]

Finally, we introduce
\[
c_{j,i}=w_{j,i}-a_i-\frac{b_i}j=-b_i\cdot\frac{\beta+1}{j(\alpha
j+\beta+1)}\,.
\]

Since $(a_i)$ is a probability distribution, for (\ref{iii}) it suffices
to show that $(a_i)$, $(b_i)$, $(c_{j,i})$, and $(r_i)$ are
exponentially small. 

According to Assumption \ref{puba7}, $Y$ has finite moment 
generating function. This implies  that $\lef a_i\ri$ is exponentially
small. The same holds for $\lef b_i\ri$, because
\[
\sum_{i=1}^\infty H(i)e^{\veps i}=
\E\bigl((\nu-1)e^{\veps Y})\bigr)\le\Bigl[\E(\nu-1)^2\;
\E\bigl(e^{2\veps Y}\bigr)\Bigr]^{\!1/2}<\infty
\]
if $\veps>0$ is small enough. Finally,
\begin{multline*}
\sum_{j=1}^\infty\sum_{i=1}^{j-1}|c_{j,i}|e^{\veps i}=
\sum_{j=1}^\infty\sum_{i=1}^{j-1}|b_i|\,\frac{\beta+1}{j(\alpha
j+\beta+1)}\,e^{\veps i}\\
\le\sum_{j=1}^\infty \frac{\beta+1}{j(\alpha j+\beta+1)}\;
\sum_{i=1}^\infty |b_i|e^{\veps i}<\infty.
\end{multline*}

The sequence $\lef r_j\ri$ is also exponentially small, because $X$
has finite moment generating function by Assumption \ref{puba7}. 

$w_{j,i}$, $a_j$, $r_j$ are nonnegative. Assumption \ref{puba8}
guarantees that the greatest common divisor of the set $\lk j:
a_j>0\rk$ is equal to 1, and $r_j>0$ for some $j$.  

We have checked all conditions of Theorem 1 of \cite{rek}. Since $X$
is not identically $0$, there exists a $k$ with $x_k>0$. On the other
hand, by Assumption \ref{puba8}, $P\lef Y=\ell\ri>0$ for some
$\ell$. Now, one can see from the recursion that $x_k, x_{k+l},
x_{k+2l},\ldots$ are all 
positive, hence the  sequence $\lef x_n\ri$ has infinitely many positive terms. 
Therefore, applying the theorem we obtain that $x_j\sim
C\;j^{-\gamma}$ as $j\tv$, where 
\[
\gamma=-\frac{\sum_{i=1}^\infty b_i}{\sum_{i=1}^\infty ia_i}\,.
\]
It is easy to see that 
\begin{gather*}
\sum_{i=1}^\infty ia_i=\sum_{i=1}^{\infty} i\,\Prob\lef Y=i\mid
Y>0\ri=\frac{\E Y}{\Prob(Y>0)}\,;\\  
-\sum_{i=1}^\infty b_i=-\frac{\beta}{\alpha}+\sum_{i=1}^\infty ia_i+
\frac{\beta+1}{\alpha}=\frac{\E X+\E Z+\E Y}{\Prob(Y>0)}\,.
\end{gather*}
Hence the statement of Theorem \ref{pubt2} follows.
\qed

\subsection*{Proof of Theorem \ref{pubt3}}

We will use the results of the discrete part, namely, Theorem \ref{pubt1}.
Let $h$ be sufficiently small positive number. We will consider limits
as $h\rightarrow 0$.  

Let $F_n(t)=\Prob(Y_n>t)$ and $H_n(t)=\E\bigl((\nu_n-1)\mathbb
I(Y_n>t)\bigl)$, as before. Furthermore, for a decreasing function
$\vph$ let $\Delta_h\vph(t)=\vph(t-h)-\vph(t)$.
 
By Lemma \ref{publ2}, the conditional probability of the event that 
an author of weight between $t-ih$ and $t-(i-1)h$ is chosen, and his
bonus is at least $(i-1)h$, given $\mathcal F_{n-1}^+$, is bounded
from above by 
\begin{equation}\label{pube2}
\biggl[\frac{t-(i-1)h}{S_{n-1}}+
\biggl(1-\frac{t-(i-1)h}{S_{n-1}}\biggr)\frac{\nu_n-1}{n-1}\biggr]
\Prob\bigl(Y_n>(i-1)h\bigm|\Cal F_{n-1}^+\bigr).
\end{equation}
Hence the conditional probability with respect to $\mathcal F_{n-1}$
is at most
\begin{equation*}
u_i:=\frac{t-(i-1)h}{S_{n-1}}\;F_n\bigl((i-1)h\bigr)
+\frac{1}{n-1}\biggl(1-\frac{t-(i-1)h}{S_{n-1}}\biggr)
H_n\bigl((i-1)h\bigr).
\end{equation*}
Note that $u_i$ depends on $n$, which is fixed at the moment.
We get that
\begin{multline*}
\E\bigl(\xi_n(t)\bigm|\mathcal F_{n-1}\bigr)\\
\le\xi_{n-1}(t)+\sum_{i=1}^{\lceil t/h\rceil}
\Bigl[\xi_{n-1}(t-ih)-\xi_{n-1}\bigl(t-(i-1)h\bigr)\Bigr]u_i
+\Prob(X>t).
\end{multline*}
After rearranging we obtain that
\begin{multline}\label{pube4}
\E\bigl(\xi_n(t)\bigm|\mathcal F_{n-1}\bigr)\le\xi_{n-1}(t)(1-u_1)\\
+\sum_{i=1}^{\lceil t/h\rceil}\xi_{n-1}(t-ih)(u_i-u_{i+1})
+nu_{\lceil t/h\rceil+1}+\Prob(X>t).
\end{multline}

Here
\begin{multline*}
u_1=\frac{t}{S_{n-1}}+\frac{1}{n-1}\biggl(1-\frac{t}{S_{n-1}}\biggr)
\E(\nu_n-1)\\
=\biggl(\frac{t}{\E X+\E Z}+\E\nu-1\biggr)\frac{1+o(1)}{n}\,,
\end{multline*}
and
\begin{multline*}
u_i-u_{i+1}=\frac{h}{S_{n-1}}\,F_n\bigl((i-1)h\bigr)+
\frac{t-ih}{S_{n-1}}\,\Delta_hF_n(ih)\\
-\frac{1}{n-1}\;\frac{h}{S_{n-1}}\,H_n\bigl((i-1)h\bigr)
+\frac{1}{n-1}\biggl(1-\frac{t-ih}{S_{n-1}}\biggr)
\Delta_hH_n(ih).
\end{multline*}

This implies that
\[
n(u_i-u_{i+1})\to\frac{h}{\E X+\E Z}\,F\bigl((i-1)h\bigr)
+\frac{t-ih}{\E X+\E Z}\,\Delta_hF(ih)+\Delta_hH(ih),
\]
as $n\to\infty$. Finally, 
\[
nu_{\lceil t/h\rceil+1}\le\frac{n}{n-1}\Bigl(1+\frac{h}{S_{n-1}}\Bigr)
H_n(t),
\]
hence
\[
\limsup_{n\to\infty}nu_{\lceil t/h\rceil+1}\le H(t).
\]
Let
\[
G_u(t)=\limsup_{n\to\infty}\frac{\xi_n(t)}{n}\,
\]
(subscript $u$ stands for ``upper'').
$G_u(t)$ is a decreasing random function, and
\begin{multline*}
\limsup_{n\to\infty}\sum_{i=1}^{\lceil t/h\rceil}
\xi_{n-1}(t-ih)(u_i-u_{i+1})\\
\le \sum_{i=1}^{\lceil t/h\rceil}G_u(t-ih)\biggl[
\frac{F\bigl((i-1)h\bigr)}{\E X+\E Z}\,h
+\frac{t-ih}{\E X+\E Z}\,\Delta_hF(ih)+\Delta_hH(ih)\biggl].
\end{multline*}

Denote the sum on the right hand side by $\Sigma_u(t,h)$. 
We want to apply Lemma \ref{publ1} to the sequence $\xi_n(t)$.
It satisfies \eqref{pube4}, and, similarly to the discrete case,
\[
\E\bigl((\xi_n(t)-\xi_{n-1}(t))^2\bigm|\Cal
F_{n-1}\bigr)\le \E(\nu_n+1)^2=O(1)
\]
holds again. The other assumptions are also easy to check. 
Hence
\[
G_u(t)\le\Bigl[\Sigma_u(t,h)+H(t)+\Prob(X>t)\Bigr]\;\biggl[
\frac{t}{\E X+\E Z}+\E\nu\biggr]^{\!-1}.
\]
One can readily verify that $\Sigma_u(t,h)$ converges to 
\begin{multline*}
\frac{1}{\E X+\E Z}\Biggl[\int_0^t G_u(t-s)F(s)\,ds
-\int_0^t G_u(t-s)(t-s)\,dF(s)\Biggr]\\
-\int_0^t G_u(t-s)\,dH(s)=\int_0^t G_u(t-s)\,d_sL(t,s)
\end{multline*}
as $h\to 0$, since the Riemann--Stieltjes integrals in the expression
exist. This implies that
\begin{equation}
G_u(t)\le\biggl[\int_0^t G_u(t-s)\,d_sL(t,s)+H(t)+\Prob(X>t)\biggr]
\biggl[\frac{t}{\E X+\E Z}+\E\nu\biggr]^{\!-1}.\label{pube5}
\end{equation}

Therefore the solution of the corresponding integral equation
\eqref{pube0} with initial condition $G_u(0)=1$ is an upper bound for
$G_u(t)$. That is, $G_u(t)\le G(t)$, where $G(t)$ is  the
deterministic function given in the theorem.  
 
Now we give lower bounds by analogous argumentation.

We estimate from below the conditional probability that an author 
with weight between $t-ih$ and $t-(i-1)h$ is chosen and his bonus is
at least $ih$, given $\Cal F_{n-1}^+$. Similarly to
\eqref{pube2}, we have that it is greater than or equal to 
\[
\biggl[\frac{t-ih}{S_{n-1}}+
\biggl(1-\frac{t-ih}{S_{n-1}}\biggr)\frac{\nu_n-1}{n-1}\biggr]
\Prob\bigl(Y_n>ih\bigm|\Cal F_{n-1}^+\bigr).
\]
Hence the lower bound of the conditional probability with respect to
$\Cal F_{n-1}$ is the following.
\[
\ell_i:=\frac{t-ih}{S_{n-1}}\;F_n(ih)
+\frac{1}{n-1}\biggl(1-\frac{t-ih}{S_{n-1}}\biggr)H_n(ih).
\]
We obtain that
\begin{multline*}
\E\bigl(\xi_n(t)\bigm|\Cal F_{n-1}\bigr)\\ 
\ge\xi_{n-1}(t)+\sum_{i=1}^{\lceil t/h\rceil}
\Bigl[\xi_{n-1}(t-ih)-\xi_{n-1}\bigl(t-(i-1)h\bigr)\Bigr]\ell_i
+\Prob(X>t).
\end{multline*}
After rearranging we get a formula similar to \eqref{pube4}. 
\begin{multline}\label{pube6}
E\bigl(\xi_n(t)\bigm|\Cal F_{n-1}\bigr)\ge\xi_{n-1}(t)(1-\ell_1)\\
+\sum_{i=1}^{\lceil t/h\rceil}\xi_{n-1}(t-ih)(\ell_i-\ell_{i+1})
+n\ell_{\lceil t/h\rceil+1}+\Prob(X>t).
\end{multline}

Here
\begin{multline*}
\ell_1=\frac{t-h}{S_{n-1}}\,F_n(h)+\frac{1}{n-1}
\biggl(1-\frac{t-h}{S_{n-1}}\biggr)H_n(h)\\
=\biggl(\frac{t-h}{\E X+\E Z}\,F(h)+H(h)\biggr)\frac{1+o(1)}{n}\,,
\end{multline*}
and
\begin{multline*}
\ell_i-\ell_{i+1}=\frac{h}{S_{n-1}}\,F_n(ih)+
\frac{t-(i+1)h}{S_{n-1}}\,\Delta_hF_n\bigl((i+1)h\bigr)\\
-\frac{1}{n-1}\;\frac{h}{S_{n-1}}\,H_n(ih)
+\frac{1}{n-1}\biggl(1-\frac{t-(i+1)h}{S_{n-1}}\biggr)
\Delta_hH_n\bigl((i+1)h\bigr).
\end{multline*}
This implies that $n(\ell_i-\ell_{i+1})$ converges to
\[
\frac{h}{\E X+\E Z}\,F(ih)
+\frac{t-(i+1)h}{\E X+\E Z}\,\Delta_hF\bigl((i+1)h\bigr)+
\Delta_hH\bigl((i+1)h\bigr)
\]
as $n\to\infty$. Finally,
\[
n\ell_{\lceil t/h\rceil+1}\ge-\frac{2nh}{S_{n-1}}+
H_n(t+2h),
\]
therefore
\[
\liminf_{n\to\infty}n\ell_{\lceil t/h\rceil+1}\ge-
\frac{2h}{\E X+\E Z}+H(t+2h).
\]
Let
\[
G_\ell(t)=\liminf_{n\to\infty}\frac{\xi_n(t)}{n}\,;
\]
then $G_\ell(t)$ is also a decreasing random function. On the right
hand side of \eqref{pube6} we have
\[
\liminf_{n\to\infty}\sum_{i=1}^{\lceil t/h\rceil}
\xi_{n-1}(t-ih)(\ell_i-\ell_{i+1})\ge \Sigma_\ell(t,h),
\]
where
\begin{multline*}
\Sigma_\ell(t,h)=\sum_{i=1}^{\lceil t/h\rceil}G_\ell(t-ih)\biggl[
\frac{F(ih)}{\E X+\E Z}\,h\\
+\frac{t-(i+1)h}{\E X+\E Z}\,\Delta_hF\bigl((i+1)h\bigr)
+\Delta_hH\bigl((i+1)h\bigr)\biggl].
\end{multline*}

Applying Lemma \ref{publ1} we get that
\begin{multline*}
G_\ell(t)\ge \biggl[\Sigma_\ell(t,h)-\frac{2h}{\E X+\E Z}+H(t+2h)+\Prob(X>t)
\biggr]\\
\times\biggl[\frac{t-h}{\E X+\E Z}\,F(h)+H(h)+1\biggr]^{\!-1}.
\end{multline*}

Let  $h$ go to  zero again. The sum $\Sigma_\ell(t,h)$
converges to the same Riemann--Stieltjes integral as $\Sigma_u(t,h)$
does. Thus the right hand side of the inequality above converges to 
the right hand side of \eqref{pube5}. Hence we obtain that
$G_\ell(t)\ge G(t)$. This, together with the estimation for $G_u(t)$,
implies the statement of the theorem. 
\kesz

\subsection*{Proof of Theorem \ref{pubt4}}


Let the density function of $Y$ be denoted by $f$. From the absolute
continuity of $F$ the same follows for $H$. Let $h$ be defined by 
\[
H(t)=\int_t^\infty h(s)\,ds.
\]
Differentiating $L$
with respect to $s$ we obtain that  
\[
\frac{\partial}{\partial s}L\lef t, s\ri=\frac{F\lef s \ri-sf\lef
s\ri+tf\lef s\ri}{\E X+\E Z}+h\lef s\ri \quad \lef0\le s\le
t\ri\,.
\] 

Hence equation \eqref{pube0} may be written in the following form. 
\[G\lef t\ri=\int_0^tG\lef t-s\ri w_{t,s}ds+r\lef t\ri,\]
where
\begin{align*}
w_{t,s}&=\frac{\dfrac{F\lef s\ri+\lef t-s\ri f\lef s\ri}
{\E X+\E Z}+h(s)}{\dfrac{t}{\E X+\E Z}+\E \nu}\\
&=\frac{F\lef s\ri+\lef t-s\ri f\lef s\ri+h\lef s\ri\lef\E X+\E Z\ri}
{t+\lef \E X+\E Z\ri\E \nu}\,;\\ 
r\lef t\ri&=\frac{H\lef t\ri+\Prob\lef X>t\ri}
{\dfrac{t}{\E X+\E Z}+\E \nu}.
\end{align*}

In order to apply Theorem 2 of \cite{rek} write $w_{t,s}$ in the
following form.
\begin{align*}
w_{t,s}&=f\lef s\ri+\frac{F\lef s\ri-\lef s+\lef \E X+\E Z\ri
\E \nu \ri f\lef s\ri+h\lef s\ri\lef \E X+\E Z\ri}
{t+\lef \E X+\E Z\ri\E \nu}\\
&=f\lef s\ri+\frac{b\lef s\ri}{t+d}\,,
\end{align*}
where
\begin{gather*}
b\lef s\ri=F\lef s\ri-\bigl(s+\lef \E X+\E Z\ri\E \nu \bigr)
f\lef s\ri+h\lef s\ri\lef \E X+\E Z\ri;\\ 
d=\lef \E X+\E Z\ri\E \nu.
\end{gather*}

Next we check that all assumptions required in \cite{rek} hold. Since 
$f$ is a probability density function, $G$ is clearly decreasing and
$w$ is nonnegative, all we need is the following three facts.
\begin{enumerate}[(i)]
\item\label{2-i}
$d$ is a positive constant,
\item
$r$ is a nonnegative, continuous function,
\item\label{2-iii}
there exists $z>1$ such that 
\begin{gather*}
\int_0^{\infty} f\lef t\ri z^t dt<\infty, \qquad
\int_0^{\infty} \left\vert b\lef t\ri \right\vert z^t dt<\infty,
\end{gather*}
and $r\lef t\ri z^t$ is directly Riemann integrable on $\lsz
0,\infty\ri$.
\end{enumerate}

Here (\ref{2-i}) follows from Assumption \ref{puba8}.  
From the continuity of $F$ and $H$ the same follows for $r$.
Finally, the first part of condition (\ref{2-iii}) easily follows from
Assumptions 2 and 7. In addition, using that $r$ is monotonically
decreasing we get that  
\[
\sum_{n=1}^{\infty}\sup_{0\leq \theta\leq \tau} r\lef
t+n\tau+\theta\ri z^{t+n\tau+\theta}\leq \sum_{n=1}^{\infty}\lsz r\lef
t+n\tau\ri z^{t+n\tau}\rsz z^{\tau}
\] 
for $z>1$. The right hand side is finite for almost all $t$, because
$\int_0^{\infty}r\lef s\ri z^s ds$ is finite. Therefore $r\lef t\ri
z^t$ is directly Riemann integrable.  

Thus Theorem \ref{pubt4} follows from Theorem 2 of \cite{rek}.
Using the continuity of $G$ and the method of the discrete case it is
easy to see that $G$ is not identically $0$ for large $t$, thus it is
polynomially decaying. What is left is to determine the exponent, that
is,  
\[
\gamma=-\frac{\int_0^{\infty}b\lef s\ri ds}{\int_0^{\infty} s f
\lef s\ri ds}.
\] 
The denominator is equal to $\E Y$. In the numerator we have 
\begin{align*}
&\int_0^{\infty}b\lef s\ri ds\\
&=\int_0^{\infty}\Bigl(F\lef s\ri-\bigl(s+\lef \E X+\E Z\ri\E\nu\bigr)
f\lef s\ri+
h\lef s\ri\lef \E X+\E Z\ri\Bigr)ds\\
&=\E Y-\E Y- \lef \E X+\E Z\ri\E \nu +H(0)\lef \E X+\E Z\ri\\ 
&=- \lef \E X+\E Z\ri\E \nu+\E\lef \nu-1\ri \lef \E X+\E Z\ri\\
&=-\lef\E X+\E Z\ri. 
\end{align*}

Therefore we got that 
\[
\gamma=\frac{\E X+\E Z}{\E Y},
\]
and the proof of Theorem \ref{pubt4}  is complete. \qed


\begin{thebibliography}{99}
\bibitem{rek} \'A.\ Backhausz, T.\ F.\ M\'ori, Asymptotics of a 
renewal-like recursion and an integral equation, {\it Appl. Anal. Discrete Math.} {\bf 8} (2014), 200--223.

\bibitem {ba} A-L.\ Barab{\'a}si, R.\ Albert, Emergence
of scaling in random networks, \textit{Science} \textbf{286} (1999), 
509--512.  

\bibitem{bingham} N.\ H.\ Bingham, C.\ M.\ Goldie, J.\ L.\ Teugels,
\textit{Regular variation}, Encyclopedia of Mathematics and its
Applications, 27, Cambridge Univ.\ Press, Cambridge, 1987. MR0898871
  (88i:26004) 

\bibitem{bojanic} R.\ Bojani\'c, E.\ Seneta, Slowly varying
functions and asymptotic relations, \textit{J.\ Math.\ Anal.\ Appl.}, 
{\bf 34} (1971), 302--315. MR0274676 (43 \#438)  
\bibitem{chung} F.\ Chung, L.\ Lu, {\it Complex graphs and
networks}, CBMS Regional Conference Series in Mathematics, 107,
Published for the Conference Board of the Mathematical Sciences,
Washington, DC, 2006. MR2248695 (2007i:05169) 

\bibitem{[Ne75]} J.\ Neveu, \textit{Discrete-parameter martingales.}
North-Holland Publishing Co., New York, 1975. MR0402915  

\bibitem{petrov} V.\ V.\ Petrov, \textit{Limit theorems of probability
theory}, Oxford Univ.\ Press, New York, 1995. MR1353441 (96h:60048) 

\bibitem{pittel} B.\ Pittel, Note on the heights of random
recursive trees and random $m$-ary search trees, \textit{Random
Struct.\ Algorithms} \textbf{5} (1994), 337--348.



\end{thebibliography}
\end{document}